\documentclass[10pt]{article}

 \newtheorem{theorem}{Theorem}

 \newtheorem{definition}{Definition}
 \newtheorem{algorithm}{Algorithm}
\begin{document}

\title{Multivariate Polynomial Factorization by Interpolation Method
\thanks{The work is partially supported by China 973 Project NKBRPC-2004CB318003.}}
\author{
 Jingzhong Zhang\ \ \ \ Yong Feng \\ \\Xijin Tang \\
Laboratory for Automated Reasoning and Programming\\
Chengdu Institute of Computer Applications\\
Chinese Academy of Sciences\\
610041 Chengdu, P. R. China\\
E-mail: zjz101@yahoo.com.cn, yongfeng@casit.ac.cn }
\date{}
\maketitle

\begin{abstract}
Factorization of polynomials arises in numerous areas in symbolic
computation. It is an important capability in many symbolic and
algebraic computation. There are two type of factorization of
polynomials. One is convention polynomial factorization, and the
other approximate polynomial factorization.

Conventional factorization algorithms use symbolic methods to get
exact factors of a polynomial while approximate factorization
algorithms use numerical methods to get approximate factors of a
polynomial. Symbolic computation often confront intermediate
expression swell problem, which lower the efficiency of
factorization. The numerical computation is famous for its high
efficiency, but it only gives approximate results. In this paper,
we present an algorithm which use approximate method to get exact
factors of a multivariate polynomial. Compared with other methods,
this method has the numerical computation advantage of high
efficiency for some class of polynomials with factors of lower
degree. The experimental results show that the method is more
efficient than {\it factor} in Maple 9.5 for polynomials with more
variables and higher degree.

\end{abstract}

\section{Introduction}
Polynomial factorization plays a significant role in many problems
including the simplification, Gr\"obner basis and solving
polynomial equations etc. It has been studied for a long time and
some high efficient algorithms have been proposed. There are two
type of factorization of polynomials. One is convention polynomial
factorization,  and the other is approximate polynomial
factorization.

The modern conventional factorization methods follow Zassenhaus'
approach \cite{Zassenhaus_a}\cite{Zassenhauss_b}. First,
Multivariate polynomial factorization is reduced to bivariate
factorization due to Bertini's theorem and hensel
lifting\cite{GAO}\cite{GATHEN_a}. Then one of the two remaining
variables is specialized at random. The resulting univariate
polynomial is factored and its factors are lifted up to a high
enough precision. At last, the lifted factors are recombined to
get the factors of the original polynomial.

Approximate factorization is a natural extension of conventional
polynomial factorization. It uses approximate methods to get
approximate factorization of polynomial. The approximate
factorization is not popular now, but there are some papers to
discuss it. In 1985, Kaltonfen presented an algorithm for
performing the absolute irreducible factorization, and suggested
to perform his algorithm by floating-point numbers, then the
factor obtained is an approximate one. However, the concept of
approximate factorization appeared first in a paper on control
theory\cite{Mou}. The algorithm is as follows: At first express
the two factors $G$ and $H$ of the polynomials $F$ with unknown
coefficients by fixing their terms, then determine the numerical
coefficients so as to minimize $\|F-GH \|$. Huang et al. pursuit
this approach, but the algorithm seems to be rarely successful,
unless $G$ or $H$ is a polynomial of several terms. In 1991,
Sasaki et al. proposed a modern algorithm\cite{Sasaki_a}, which
use power-series roots to find approximate factors. This algorithm
is successful for polynomials of small degrees. Subsequently,
Sasaki et al. presented another algorithm\cite{Sasaki_b} which
utilizes zero-sum relations. The zero-sum relations are quite
effective for determining approximate factors. However,
computation based on zero-sum relations is practically very
time-consuming. In \cite{Sasaki_c}, Sasaki, T. presented an
effective method to get as many zero-sum relations as possible by
matrix operations so that approximate factorization algorithm is
improved. In \cite{corless_B}, Corless et al proposed an algorithm
for factoring bivariate approximate polynomial based on the idea
of decomposition of affine variety. However, it is not easy to
generalize the algorithm to factor multivariate approximate
polynomials. Recently, Zhang et al\cite{zhang} proposed an
algorithm for obtaining exact rational number from its approximate
floating number. In this paper, basing on the algorithm in
\cite{zhang}, we present an algorithm which use approximate method
to get exact factors of a polynomial. It can be regard as a
generalization of Corless' algorithm in multivariate polynomial
and exact polynomial case.

The remainder of the paper is organized as follows. Section 2
presents a formula by which a polynomial is constructed from
sampled points on its variety. A condition is given for the
formula to determine only one polynomial up to a nonzero constant
multiples, and the error estimation is discussed. Section 3 gives
a review of  a modified continued fraction method, by which an
exact rational number can be obtained from its approximation.
Section 4 first discusses the error control, and then proposes a
factorization algorithm for multivariate polynomial over rational
number field. Section 5 gives some experimental results. The final
section makes conclusion.

\section{Interpolation method}
Polynomial interpolation is a classical numerical method. It is
studied very well for univariate polynomials. In general, there
are four types of polynomial interpolation method: Lagrange
Interpolation, Neville's Interpolation, Newton's Interpolation and
Hermite Interpolation. Lagrange interpolation formula can get the
interpolation polynomial at once for a given set of distinct
interpolation points and corresponding values $(x_i,f_i),
i=1,\cdots,n+1$. It is very useful in some situations in which
many interpolation problems are to be solved for the same set of
interpolation points $x_i,i=1,\cdots,n+1$, but different sets of
function values $f_i,i=1,\cdots,n+1$. Unlike Lagrange
interpolation method which solve the interpolation problem all at
once, Neville's interpolation method solve the problem for smaller
sets of interpolation points first and then update these solutions
to obtain the solution to the full interpolation problem. It aims
at determining the value of the interpolating polynomial at some
point. It is less suited for determining the interpolating
polynomial. If the interpolating polynomial is needed, Newton's
Interpolation formula is preferred. Just like Neville's
interpolation method, it first get interpolating polynomial for
smaller sets of interpolation points and then update the
polynomial for a larger sets of interpolation points, step by
step, and finally, the interpolating polynomial is obtained for
the set of the full interpolation points. If the interpolating
problem prescribes at each interpolation point
$x_i,i=1,\cdots,n+1$ not only the value but also the derivatives
of desired polynomial, then the Hermite formula is preferred.

For univariate polynomial interpolation, $n+1$ distinct
interpolation points and their corresponding values determine only
one polynomial with degree less than or equal to $n$. However, the
interpolation points for multivariate polynomial interpolation can
not be chosen arbitrarily. They need to satisfy some conditions.
So we need a definition as follows:
\begin{definition}\label{def:1}
Let $\Theta$ be a set of n-dimension points and $P$ a polynomial
space. We call $\Theta$ {\bf Proper interpolation points} of $P$
if for any $f$ defined on $\Theta$, there is a unique polynomial
$p\in P$ matching $f$ at $\Theta$.
\end{definition}

In definition \ref{def:1}, a polynomial $p$ matching $f$ at
$\Theta$ means $p|\Theta =f|\Theta$. In general, we can determine
interpolation polynomial space such as: if knowing the total
degree $d$ of $f$, we choose $P=\{p|deg(p)\le d\ ,p\in
K[x_1,x_2,\cdots,x_n] \}$; if knowing the degree $d_i$ of $f$ in
$x_i(i=1,\cdots,n)$, we choose
$P=\{p|\wedge_{i=1,2,\cdots,n}(deg_{x_i}(p)\le d_i), p\in
K[x_1,x_2,\cdots,x_n]\}$. Once the interpolation polynomial space
is determined, the proper interpolation points $\Theta$ of $P$ can
be set by interpolation methods\cite{Kaltofen}\cite{Boor}.

In this paper, we need to construct a polynomial from some points
of its variety. Values of the polynomial at interpolation points
are all zero. So, we introduce an interpolation formula for this
case.

Let $f(x_1,x_2,\cdots,x_n)$ be a polynomial to be interpolated. It
is represented as follows:
\begin{equation}\label{equ:add2}
 f(x_1,x_2,\cdots,x_n)=c_1X^{\alpha_1}+c_2X^{\alpha_2}+\cdots+c_mX^{\alpha_m},
\end{equation}
where $X^{\alpha_i}=x_1^{d_{1,i}}x_2^{d_{2,i}}\cdots
x_n^{d_{n,i}}$ are the distinct monomials, and $c_i$ are the
corresponding coefficients.

Let $p_1,p_2,\cdots,p_{m-1}$ be points on variety of
$f(x_1,x_2,\cdots,x_n)$, where $p_i=(p_{11},p_{12},\cdots,p_{1n})$
for $i=1,\cdots,m-1$.
$P_i^{\alpha_j}=p_{i1}^{d_{1,j}}p_{i2}^{d_{2,j}}\cdots
p_{in}^{d_{n,j}}$ denote the value of the monomial $X^{\alpha_j}$
at $p_i$. An interpolation formula is as follows:
\begin{equation}\label{equ:add1}
G(x_1,x_2,\cdots,x_n)=\left|
\begin{array}{llll}
X^{\alpha_1} & X^{\alpha_2} & \cdots & X^{\alpha_m} \\
P_1^{\alpha_1}& P_1^{\alpha_2}&\cdots &P_1^{\alpha_m} \\
\cdots&\cdots&\cdots&\cdots \\
P_{m-1}^{\alpha_1}&P_{m-1}^{\alpha_2}&\cdots&P_{m-1}^{\alpha_m}
\end{array}
\right|
\end{equation}

Next, we need to know what condition the interpolation points
$p_1,p_2,\cdots,p_{m-1}$ should satisfy so as to ensure formula
(\ref{equ:add1}) to determine a unique polynomial and it is
$f(x_1,x_2,\cdots)$ up to a nonzero constant multiplies.

Let $X_i^*$ denote the minor of $X^{\alpha_i}$ resulting from the
deletion of row 1 and column $i$ in formula (\ref{equ:add1}). We
have a theorem as follows.
\begin{theorem}\label{theo:1}
Let $f(x_1,x_2,\cdots,x_n)$ be a nonzero polynomial and it is
expressed as in equation (\ref{equ:add2}). If the $m-1$ zeroes of
$f(x_1,x_2,\cdots,x_n)$ satisfy a condition that
$G(x_1,\cdots,x_n)\neq 0$ in formula \ref{equ:add1}, then formula
(\ref{equ:add1})determines a unique polynomial and it is
polynomial $f(x_1,x_2,\cdots,x_n)$ up to a nonzero constant
multiples.
\end{theorem}
{\bf Proof}: Due to $f(x_1,x_2,\cdots,x_n)\neq 0$, the
coefficients of $f(x_1,x_2,\cdots,x_n)$ are not all zero. Assume
that $c_{i_1}\neq 0, c_{i_2}\neq 0,\cdots, c_{i_s}\neq 0 $, and
their corresponding monomials are
$X^{\alpha_{i_1}},X^{\alpha_{i_2}},\cdots, X^{\alpha_{i_s}}$. Let
$X_{i_1}^*,X_{i_2}^*,\cdots, X_{i_s}^*$ denote the minors of
$X^{\alpha_{i_1}},X^{\alpha_{i_2}},\cdots, X^{\alpha_{i_s}}$ in
formula (\ref{equ:add1}) respectively.

First, we show that if one of $X^*_{i1},X^*_{i2},\cdots, X^*_{is}$
is nonzero, then
$f(x_1,x_2,\cdots,x_n)=cG(x_1,x_2,\cdots,x_n)$,where $c$ is a
nonzero constant. Without loss of generality, let us assume that
$c_k\neq 0$ and $X^*_k\neq 0$. So we have
\begin{eqnarray*} &&G(x_1,x_2,\cdots,x_n)\\ &=&\left|
\begin{array}{llllllll}
X^{\alpha_1} & X^{\alpha_2} & \cdots& X^{\alpha_{k-1}} &X^{\alpha_k} &X^{\alpha_{k+1}} &\cdots & X^{\alpha_m} \\
P_1^{\alpha_1}& P_1^{\alpha_2}&\cdots&P_1^{\alpha_{k-1}}& P_1^{\alpha_k}&P_1^{\alpha_{k+1}}&\cdots &P_1^{\alpha_m} \\
\cdots&\cdots&\cdots&\cdots&\cdots&\cdots&\cdots&\cdots \\
P_{m-1}^{\alpha_1}&P_{m-1}^{\alpha_2}&\cdots&P_{m-1}^{\alpha_{k-1}}&P_{m-1}^{\alpha_k}&P_{m-1}^{\alpha_{k+1}}&\cdots&P_{m-1}^{\alpha_m}
\end{array}
\right|\\
&=&\frac{1}{c_k}\left|
\begin{array}{llllllll}
X^{\alpha_1} & X^{\alpha_2} & \cdots& X^{\alpha_{k-1}} &c_kX^{\alpha_k} &X^{\alpha_{k+1}} &\cdots & X^{\alpha_m} \\
P_1^{\alpha_1}& P_1^{\alpha_2}&\cdots&P_1^{\alpha_{k-1}}& c_kP_1^{\alpha_k}&P_1^{\alpha_{k+1}}&\cdots &P_1^{\alpha_m} \\
\cdots&\cdots&\cdots&\cdots&\cdots&\cdots&\cdots&\cdots \\
P_{m-1}^{\alpha_1}&P_{m-1}^{\alpha_2}&\cdots&P_{m-1}^{\alpha_{k-1}}&c_kP_{m-1}^{\alpha_k}&P_{m-1}^{\alpha_{k+1}}&\cdots&P_{m-1}^{\alpha_m}
\end{array}
\right|
\end{eqnarray*}
For $i=1,2,\cdots,k-1,k+1,\cdots,m$, addition of $c_i$ times
column $i$ to column $k$ yields:
\begin{eqnarray*} &&G(x_1,x_2,\cdots,x_n)\\ &=&\frac{1}{c_k}\left|
\begin{array}{lllllll}
X^{\alpha_1} & \cdots& X^{\alpha_{k-1}} &\sum_{i=1}^{k}c_iX^{\alpha_k} &X^{\alpha_{k+1}} &\cdots & X^{\alpha_m} \\
P_1^{\alpha_1}&\cdots&P_1^{\alpha_{k-1}}& \sum_{i=1}^{k}c_iP_1^{\alpha_k}&P_1^{\alpha_{k+1}}&\cdots &P_1^{\alpha_m} \\
\cdots&\cdots&\cdots&\cdots&\cdots&\cdots&\cdots \\
P_{m-1}^{\alpha_1}&\cdots&P_{m-1}^{\alpha_{k-1}}&\sum_{i=1}^{k}c_iP_{m-1}^{\alpha_k}&P_{m-1}^{\alpha_{k+1}}&\cdots&P_{m-1}^{\alpha_m}
\end{array}
\right|\\
&=&\frac{1}{c_k}\left|
\begin{array}{cccccccc}
X^{\alpha_1} & \cdots& X^{\alpha_{k-1}} &f(x_1,\cdots,x_n) &X^{\alpha_{k+1}} &\cdots & X^{\alpha_m} \\
P_1^{\alpha_1}&\cdots&P_1^{\alpha_{k-1}}& 0&P_1^{\alpha_{k+1}}&\cdots &P_1^{\alpha_m} \\
\cdots&\cdots&\cdots&\cdots&\cdots&\cdots&\cdots \\
P_{m-1}^{\alpha_1}&\cdots&P_{m-1}^{\alpha_{k-1}}&0&P_{m-1}^{\alpha_{k+1}}&\cdots&P_{m-1}^{\alpha_m}
\end{array}
\right|
\\
&=&\frac{(-1)^{k+1}}{c_k}X_k^*f(x_1,\cdots,x_n)
\end{eqnarray*}

Due to $c_k\neq 0$ and $X_k^*\neq 0$, it follows that
$f(x_1,x_2,\cdots,x_n)=cG(x_1,x_2,\cdots,x_n)$, where
$c=\frac{(-1)^{k+1}c_k}{X_k^*}$ is nonzero constant.

Second, we assert that one of $X^*_{i1},X^*_{i2},\cdots, X^*_{id}$
must be nonzero. We prove it by contradiction. Let us assume that
$X^*_{i1}=0,X^*_{i2}=0,\cdots, X^*_{id}=0$.

Under the assumption of the theorem,
$G(x_1,x_2,\cdots,x_n)=\sum_{i=1}^mb_iX^{\alpha_i}\neq 0$ and
$f(x_1,x_2,\cdots,x_n)\neq 0$. So, not all of their coefficients
are zero. Assume that $b_{h_1}\neq 0,b_{h_2}\neq
0,\cdots,b_{h_z}\neq 0$ and $c_{i_1}\neq 0, c_{i_2}\neq 0,\cdots,
c_{i_s}\neq 0 $. Since $X^*_{i1}=0,X^*_{i2}=0,\cdots, X^*_{is}=0$,
it holds that $i_k\neq h_j$ for $k=1,\cdots, s$ and
$j=1,\cdots,z$. Hence, we have that $X^{\alpha_{i_k}}\neq
X^{\alpha_{h_j}}$ for $k=1,\cdots, s$ and $j=1,\cdots,z$. Let
\begin{eqnarray*}
H(x_1,x_2,\cdots,x_n)&=&f(x_1,x_2,\cdots,x_n)+G(x_1,x_2,\cdots,x_n)\\
&=&c_{i_1}X^{\alpha_{i_1}}+\cdots+c_{i_s}X^{\alpha_{i_s}}+b_{h_1}X^{\alpha_{h_1}}+\cdots+b_{h_z}X^{\alpha_{h_z}}\\
&\neq& 0
\end{eqnarray*}
Because of $b_{h_1}\neq 0$ and $X^*_{h_1}\neq 0$, it has been
shown above that
$$G(x_1,x_2,\cdots,x_n)=\bar{c}H(x_1,x_2,\cdots,x_n),$$ where
$\bar{c}$ is nonzero constant. Hence we deduce that
$$G(x_1,x_2,\cdots,x_n)-\bar{c}H(x_1,x_2,\cdots,x_n)\equiv 0$$.
However it is impossible because the term
$\bar{c}c_{i_1}X^{\alpha_{i_1}}\neq 0$ is not monomial of
$G(x_1,\cdots,x_n)$. Therefore we show that one of
$X^*_{i1},X^*_{i2},\cdots, X^*_{is}$ must be nonzero. The proof of
the theorem is finished.

The above theorem shows that if formula (\ref{equ:add1}) gives a
nonzero  polynomial, and the first row contains all monomials of
an interpolating polynomial. then it is the interpolating
polynomial up to nonzero constant multiples.

However, due to floating-point computation, we only get
approximate zeroes of $f$. Accordingly, we only obtain approximate
factors. In the remaining of this section, we study how the error
resulting from floating computation affects that of factors. For
simplicity, let us give a definition:
\begin{definition}
Let
$X_i^{\alpha_j}=x_{i,1}^{\alpha_{j,1}}x_{i,2}^{\alpha_{j,2}}\cdots
x_{i,n}^{\alpha_{j,n}}$, where $x_{i,1},\dots,x_{i,n}$ are complex
number and $\alpha_{j,1},\cdots,\alpha_{j,n}$ are nonnegative
integer. A generalized Vandemonder determinant is defined as
follows:
\begin{equation}
V_m=\left|
\begin{array}{llll}
X_1^{\alpha_1} & X_1^{\alpha_2} & \cdots & X_1^{\alpha_m} \\
X_2^{\alpha_1} & X_2^{\alpha_2} & \cdots & X_2^{\alpha_m} \\
\cdots&\cdots&\cdots&\cdots \\
X_m^{\alpha_1} & X_m^{\alpha_2} & \cdots & X_m^{\alpha_m}
\end{array}
\right|
\end{equation}
\end{definition}

We have an estimation of generalized vandemonder determinant as
follows:
\begin{theorem} \label{theo:3}
Let $M=\max_{i,j}\{X_i^{\alpha_j}\}$ and
$B=\max_{i,j,k}\{\|X_i^{\alpha_j}-X_i^{\alpha_k}\|,
|X_j^{\alpha_i}-X_k^{\alpha_i}\|\}$. Then for $m\geq 2$ it holds
that
$$|V_m|\leq m!M^{m-1}B$$
\end{theorem}
{\bf Proof}: We prove it by inductive method. When $m=2$, the
generalized Vandemonder determinant is
\begin{eqnarray*}
V_2&=&\left|
\begin{array}{ll}
X_1^{\alpha_1} & X_1^{\alpha_2} \\
X_2^{\alpha_1} & X_2^{\alpha_2}
\end{array}
\right|=X_1^{\alpha_1}X_2^{\alpha_2}-X_2^{\alpha_1}X_1^{\alpha_2}\\
&=&X_1^{\alpha_1}X_2^{\alpha_2}-X_1^{\alpha_1}X_1^{\alpha_2}
+X_1^{\alpha_1}X_1^{\alpha_2}-X_2^{\alpha_1}X_1^{\alpha_2} \\
&=&X_1^{\alpha_1}(X_2^{\alpha_2}-X_1^{\alpha_2})+X_1^{\alpha_2}(X_1^{\alpha_1}-X_2^{\alpha_1})
\end{eqnarray*}
So, $$|V_2|\leq
|X_1^{\alpha_1}(X_2^{\alpha_2}-X_1^{\alpha_2})|+|X_1^{\alpha_2}(X_1^{\alpha_1}-X_2^{\alpha_1})|\leq
2MB=2!M^{2-1}B
$$
Assume that  $|V_m|\leq m!M^{m-1}B$ for $m=k$. Let us show that it
holds for $m=k+1$. We expand $V_{k+1}$ by minors as follows
$$|V_{k+1}|=|\sum_{j=1}^{k+1}(-1)^{i+j}X_i^{\alpha_j}V_k^{i,j}|\leq \sum_{j=1}^{k+1}|X_i^{\alpha_j}|*|V_k^{i,j}|$$
According to our assumption that $|V_k^{i,j}|\leq k!M^{k-1}B$, we
have
$$|V_{k+1}|\leq \sum_{j=1}^{k+1}|X_i^{\alpha_j}|*|V_k^{i,j}|\leq \sum_{j=1}^{k+1}M*k!M^{k-1}B=(k+1)!M^kB $$
The proof is finished.

\begin{theorem}\label{theo:4}
Let $M=\max_{i,j}\{X_i^{\alpha_j}\}$ and
$B=\max_{i,j,k}\{\|X_i^{\alpha_j}-X_i^{\alpha_k}\|,
|X_j^{\alpha_i}-X_k^{\alpha_i}\|\}$ and
$\varepsilon=\max_{i=1}^m|a_i|$. A determinant is as follows:
$$V_m=\left |
\begin{array}{cccccc}
X_1^{\alpha_1} & X_1^{\alpha_2} & \cdots &a_1&\cdots & X_1^{\alpha_m} \\
X_2^{\alpha_1} & X_2^{\alpha_2} & \cdots &a_2&\cdots & X_2^{\alpha_m} \\
\cdots&\cdots&\cdots&\cdots \cdots &\cdots\\
X_m^{\alpha_1} & X_m^{\alpha_2} & \cdots &a_m&\cdots&
X_m^{\alpha_m}
\end{array} \right |
$$
Then we have an estimate that $|V_m|\leq M^{m-2}m!B\varepsilon$
for $m\geq 3$.
\end{theorem}
{\bf Proof}. Expanding $V_m$ by  column  $(a_1,\cdots, a_m)^T$ and
then using theorem \ref{theo:3}, we can get the proof.

And now, we study the difference between two generalized
Vandemonder determinants.
\begin{theorem}\label{theo:5}
Let $$V_m^{(1)}=\left|
\begin{array}{llll}
X_1^{\alpha_1} & X_1^{\alpha_2} & \cdots & X_1^{\alpha_m} \\
X_2^{\alpha_1} & X_2^{\alpha_2} & \cdots & X_2^{\alpha_m} \\
\cdots&\cdots&\cdots&\cdots \\
X_m^{\alpha_1} & X_m^{\alpha_2} & \cdots & X_m^{\alpha_m}
\end{array}
\right|$$
 and
 $$V_m^{(2)}=\left|
\begin{array}{llll}
Y_1^{\alpha_1} & Y_1^{\alpha_2} & \cdots & Y_1^{\alpha_m} \\
Y_2^{\alpha_1} & Y_2^{\alpha_2} & \cdots & Y_2^{\alpha_m} \\
\cdots&\cdots&\cdots&\cdots \\
Y_m^{\alpha_1} & Y_m^{\alpha_2} & \cdots & Y_m^{\alpha_m}
\end{array}
\right|$$ and assume that
$M=\max_{i,j}\{X_i^{\alpha_j},Y_i^{\alpha_j}\}$ ,
$B=\max_{i,j,k}\{\|X_i^{\alpha_j}-X_i^{\alpha_k}\|,
|X_j^{\alpha_i}-X_k^{\alpha_i}\|,\|Y_i^{\alpha_j}-Y_i^{\alpha_k}\|,
|Y_j^{\alpha_i}-Y_k^{\alpha_i}\|\}$ and
$\varepsilon=\max_{i,j=1}^m \{\|X_i^{\alpha_j}-Y_i^{\alpha_j}\|
\}$. Then it holds for $m\ge 3$ that
$$|V_m^{(1)}-V_m^{(2)}|\leq mm!M^{m-2}B\varepsilon $$
\end{theorem}
{\bf Proof}: \begin{eqnarray*}
&&V_m^{(1)}-V_m^{(2)}=\left|\begin{array}{llll}
X_1^{\alpha_1} & X_1^{\alpha_2} & \cdots & X_1^{\alpha_m} \\
X_2^{\alpha_1} & X_2^{\alpha_2} & \cdots & X_2^{\alpha_m} \\
\cdots&\cdots&\cdots&\cdots \\
X_m^{\alpha_1} & X_m^{\alpha_2} & \cdots & X_m^{\alpha_m}
\end{array}
\right|-\left|\begin{array}{llll}
Y_1^{\alpha_1} & Y_1^{\alpha_2} & \cdots & Y_1^{\alpha_m} \\
Y_2^{\alpha_1} & Y_2^{\alpha_2} & \cdots & Y_2^{\alpha_m} \\
\cdots&\cdots&\cdots&\cdots \\
Y_m^{\alpha_1} & Y_m^{\alpha_2} & \cdots & Y_m^{\alpha_m}
\end{array}
\right| \\
&&=\left|\begin{array}{llll}
X_1^{\alpha_1} & X_1^{\alpha_2} & \cdots & X_1^{\alpha_m} \\
X_2^{\alpha_1} & X_2^{\alpha_2} & \cdots & X_2^{\alpha_m} \\
\cdots&\cdots&\cdots&\cdots \\
X_m^{\alpha_1} & X_m^{\alpha_2} & \cdots & X_m^{\alpha_m}
\end{array}
\right|-\left|\begin{array}{llll}
Y_1^{\alpha_1} & X_1^{\alpha_2} & \cdots & X_1^{\alpha_m} \\
Y_2^{\alpha_1} & X_2^{\alpha_2} & \cdots & X_2^{\alpha_m} \\
\cdots&\cdots&\cdots&\cdots \\
Y_m^{\alpha_1} & X_m^{\alpha_2} & \cdots & X_m^{\alpha_m}
\end{array}
\right|\\
&&+\left|\begin{array}{llll}
Y_1^{\alpha_1} & X_1^{\alpha_2} & \cdots & X_1^{\alpha_m} \\
Y_2^{\alpha_1} & X_2^{\alpha_2} & \cdots & X_2^{\alpha_m} \\
\cdots&\cdots&\cdots&\cdots \\
Y_m^{\alpha_1} & X_m^{\alpha_2} & \cdots & X_m^{\alpha_m}
\end{array}
\right|-\left|\begin{array}{llll}
Y_1^{\alpha_1} & Y_1^{\alpha_2} & \cdots & X_1^{\alpha_m} \\
Y_2^{\alpha_1} & Y_2^{\alpha_2} & \cdots & X_2^{\alpha_m} \\
\cdots&\cdots&\cdots&\cdots \\
Y_m^{\alpha_1} & Y_m^{\alpha_2} & \cdots & X_m^{\alpha_m}
\end{array}
\right|\\
&&+\cdots\\
&&+\left|\begin{array}{llll}
Y_1^{\alpha_1} & Y_1^{\alpha_2} & \cdots & X_1^{\alpha_m} \\
Y_2^{\alpha_1} & Y_2^{\alpha_2} & \cdots & X_2^{\alpha_m} \\
\cdots&\cdots&\cdots&\cdots \\
Y_m^{\alpha_1} & Y_m^{\alpha_2} & \cdots & X_m^{\alpha_m}
\end{array}
\right|-\left|\begin{array}{llll}
Y_1^{\alpha_1} & Y_1^{\alpha_2} & \cdots & Y_1^{\alpha_m} \\
Y_2^{\alpha_1} & Y_2^{\alpha_2} & \cdots & Y_2^{\alpha_m} \\
\cdots&\cdots&\cdots&\cdots \\
Y_m^{\alpha_1} & Y_m^{\alpha_2} & \cdots & Y_m^{\alpha_m}
\end{array}
\right|
\end{eqnarray*}
From theorem \ref{theo:4}, it holds that
$|V_m^{(1)}-V_m^{(2)}|\leq
m*M^{m-2}m!B\varepsilon=mm!M^{m-2}B\varepsilon$. The proof is
finished.

\section{Continued fraction method}
As we said above, our method is to use approximate method to get
exact factors of a multivariate polynomial over rational number
field. So we need to recover the exact coefficients of a
polynomial from its approximate coefficients. In this section, we
introduce a continued fraction method to recover exact rational
number from its approximation. As we know, a continued fraction
representation of a real number is one of the forms:
\begin{equation}\label{equ:1}
 a_0+\frac{1}{a_1+\frac{1}{a_2+\frac{1}{a_3+\cdots}}},
\end{equation}
where $a_0$ is a integer and $a_1,a_2,a_3,\cdots$ are positive
integers. One can abbreviate the above continued fraction as
$$[a_0;a_1,a_2,\cdots]$$

In order to recover exact rational number, we introduce a control
error into the conventional continued fraction method. The
continued fraction method is modified as follows.
\newcounter{num}
\begin{algorithm}\label{algor:1} Continued fraction method \\
Input: a nonnegative floating-point number $a$ and $\varepsilon_1>0$; \\
Output: a rational number $b$
\begin{list}{Step \arabic{num}:}{\usecounter{num}\setlength{\rightmargin}{\leftmargin}}
\item $i:=1$ and $x_1:=a$; \item \label{step1:2}Getting integral
part of $x_i$ and assigning it to $a_i$, assigning its remains to
$b_i$. If $b_i<\varepsilon_1$, then goto Step \ref{step1:5}; \item
$i:=i+1$; \item $x_i:=\frac{1}{b_{i-1}}$ and goto Step
\ref{step1:2}; \item \label{step1:5}  Computing expression
(\ref{equ:1}) and assigning it to $b$. \item return $b$.
\end{list}
\end{algorithm}

In \cite{zhang}, we discussed how to get  error control
$\varepsilon_1>0$. The theorem is as follows:
\begin{theorem}\label{theo:practice_alg}
Let $n_0/n_1$ be a reduced rational number and $r$ its
approximation. Assume that $n_0$,$n_1$ are  positive integers and
$L\ge\max\{n_1,2\}$. $K$ is a positive integer. The continued
fraction representations of $n_0/n_1$ and $r$ are
$[a_0,a_1,\cdots,a_N]$ and $[b_0,b_1,\cdots,b_M]$ respectively. If
$|d|=|r-n_0/n_1|<1/((2K+2)L(L-1))$, then one of the two statement
holds
\begin{itemize}
\item $a_i=b_i$ for $i=0,\cdots,N$, and $b_{N+1}\ge K$; \item
$a_i=b_i$ for $i=0,\cdots,N-1$, and $b_N=a_N-1$, $b_{N+1}=1$,
$b_{N+2}\ge K$.
\end{itemize}
\end{theorem}

From theorem \ref{theo:practice_alg},  getting exact non-negative
number $n_2/n_ 1$ from its approximation $r_0$ is summarized as
follows:
\begin{algorithm}\label{alg:2} Obtaining Exact Number\\
\begin{list}{Step \arabic{num}:}{\usecounter{num}\setlength{\rightmargin}{\leftmargin}}
\item estimating an upper bound of the denominator of $n_2/n_1$,
Denoted by $L$; \item computing $$\beta=\frac{1}{(2L+2)L(L-1))}$$
\item obtaining $r_0$ by approximate method such that
$|r_0-n_2/n_1|<\beta$; \item taking $\varepsilon_1=1/L$ in
algorithm \ref{algor:1} and calling algorithm \ref{algor:1} to get
$b$. So $n_2/n_1=b$.
\end{list}
\end{algorithm}
\section{Factoring Multivariate Polynomials by Approximate Method}

In this paper, we only discuss factorization of a multivariate
polynomial over rational number field. So its coefficients are all
rational numbers. In order to get factors of a multivariate
polynomial over rational number field, we first compute its
factors over complex field. These factors are complex coefficient
polynomials. Products of some of them must be real polynomials. We
get the products which are approximate rational coefficient
factors of the original polynomial. Finally, transforming these
real products into rational coefficient polynomials yields factors
of the original polynomial over rational number field.

 A set
$$V(f)=\{(a_1,\cdots,a_n)\in C^n: f(a_1,\cdots,a_n)=0\}$$
is called {\bf affine variety} of $f(x_1,\cdots,x_n)$. An affine
variety $V\subset C^n$ is {\bf irreducible} if whenever $V$ is
written in the form $V=V_1\cup V_2$, where $V_1$ and $V_2$ are
affine varieties, then either $V_1=V$ or $V_2=V$.

Let $f(x_1,x_2,\cdots,x_n)$ be a square free polynomial over
complex number field $C$ and $f=f_1f_2\cdots f_m$, where $f_i$ is
distinct irreducible polynomials. Then $<f>$ is a radical ideal.
It holds as follows
\begin{equation}\label{equ:2}
V(f)=V(f_1)\cup V(f_2)\cup\cdots \cup V(f_m),
\end{equation}
 where
$V(f_i)$ are irreducible affine varieties.

From equation (\ref{equ:2}), if we get a point on variety of $f$,
it must be either on one of $V(f_i)$ or on the intersection of
them. When the point is not singular point, it must be on one of
$V(f_i)$ and not on the intersection of two varieties. Theorem
\ref{theo:1} shows that if getting enough points in some variety
of $V(f_i)$ that satisfy the condition of theorem \ref{theo:1}, we
can recover the polynomial by formula (\ref{equ:add1}). Therefore,
the procedure of factorization is as follows: First get a initial
nonsingular point on one variety of $V(f_i)$. And then obtain
enough sampled points on the same variety. Third, use formula
(\ref{equ:add1}) to get a factor and finally, obtain a rational
factor.

However, due to approximate computation, we first discuss error
control, and then study factorization.

\subsection{Error control }

Let $f(x_1,\cdots,x_n)$ be a polynomial to be factored over
rational number field. According to algorithm \ref{algor:add1},
the first thing we need to do is to determine an upper bound of
absolute values of coefficient denominators of factors of
$f(x_1,x_2,\cdots,x_n)$. The following theorem is very helpful.
\begin{theorem}\label{theo:6}
Let $g(x_1,\cdots,x_n)$ be a monic polynomial over rational number
field. Its factorization over rational number field is
$g(x_1,\cdots,x_n)=\prod_{i=1}^mh_i(x_1,\cdots,x_n)$, where all
$h_i(x_1,\cdots,x_n)$ are monic polynomials. Assume that $N>0$ is
the least common multiple of denominators of coefficients of
$g(x_1,\cdots,x_n)$. Then $N$ is an upper bound of absolute values
of denominators of coefficients of $h_i(x_1,\cdots,x_n)$ for
$i=1,\cdots,m$.
\end{theorem}
{\bf Proof}: It is clear that $Ng(x_1,\cdots,x_N)$ is  a primitive
integral coefficient polynomial. Let $N_i$ be the least common
multiple of denominators of coefficients of $h_i(x_1,\cdots,x_n)$.
Hence, $h_i(x_1,\cdots,x_n)=\frac{1}{N_i}\bar
h_i(x_1,\cdots,x_n)$,where $\bar h_i(x_1,\cdots,x_n)$ is a
primitive integral coefficient polynomial. From Gauss'lemma,
$\prod_{i=1}^m\bar h_i(x_1,\cdots,x_n)$ is a primitive polynomial
over $Z[x_1,\cdots,x_n]$. On the other hand, we have that
$$Ng(x_1,\cdots,x_n)=N\prod_{i=1}^mh_i(x_1,\cdots,x_n)=\frac{N}{\prod_{i=1}^mN_i}\prod_{i=1}^m\bar h_i(x_1,\cdots,x_n)$$
Since $Ng(x_1,\cdots,x_n)$ and $\prod_{i=1}^m\bar
h_i(x_1,\cdots,x_n)$ are primitive polynomials, it holds that
$N/(\prod_{i=1}^mN_i)=\pm 1$ . Therefore, $N=\pm
\prod_{i=1}^mN_i$. The proof is finished.

Theorem \ref{theo:6} shows that the positive least common multiple
of denominators of coefficients of a monic polynomial is also an
upper bound of absolute values of denominators of coefficients of
its monic factors.

According to algorithm \ref{alg:2} and theorem \ref{theo:5}, we
calculate control error as follows:
\begin{algorithm}\label{algor:add1} Calculating control error
\begin{list}{Step \arabic{num}:}{\usecounter{num}\setlength{\rightmargin}{\leftmargin}}
\item calculating upper bound of absolute values of the
coefficient denominators of exact factors, denoted by $L$. From
theorem \ref{theo:6}, we can take the positive least common
multiple of the coefficient denominators of a monic polynomial to
be factored. \item taking $K=L+1$ and $\varepsilon_1=1/K$ in
algorithm \ref{algor:1}. \item computing
$\beta=\frac{1}{(2K+2)L(L-1))}$. \item computing the control error
$\varepsilon$ in theorem \ref{theo:5} such that
$mm!M^{m-2}B\varepsilon\le \beta$.
\end{list}
\end{algorithm}
Therefore, in order to factor polynomial $f(x_1,\cdots,x_n)$, we
first call algorithm \ref{algor:add1} to compute control error
$\varepsilon$  and control error $\varepsilon_1$ in algorithm
\ref{algor:1}.
\subsection{Initial point}

Let $f(x_1,x_2,\cdots,x_n)$ be a square-free polynomial to be
factored over rational number field. Of course, it is a
square-free polynomial over complex number field. Choosing $n-1$
floating-point numbers $x_{1,0},x_{2,0},\cdots,x_{n-1,0}$ at
random and numerically solve
$f(x_{1,0},x_{2,0},\cdots,x_{n-1,0},x_n)=0$ for variable $x_n$
within control error $\varepsilon$. Denote this solution by
$x_{n,0}$. So the zero $(x_{1,0},x_{2,0},\cdots,x_{n,0})$ of
$f(x_1,\cdots,x_n)$ must be in one of $V(f_i)$ or the intersection
of them. If $$\nabla f(x_{1,0},x_{2,0},\cdots, x_{n,x})\neq
(0,0,\cdots,0),$$ then it must be in one of $V(f_i)$ and not in
the intersection of two of them. In order to get a neighborhood
$U_0$ of $(x_{1,0},x_{2,0},\cdots,x_{n,0})$ which is only on one
variety $V(f_i)$, we require that
\begin{equation}\label{equ:3}
\left \{
\begin{array}{c}
\frac{\partial f(x_{1,0},\cdots,x_{n,0})}{\partial x_1}\neq 0 \\
\frac{\partial f(x_{1,0},\cdots,x_{n,0})}{\partial x_2}\neq 0
\\
\cdots \cdots \cdots \\
\frac{\partial f(x_{1,0},\cdots,x_{n,0})}{\partial x_{n-1}}\neq 0
\end{array}
\right.
\end{equation}

We always assume that $\frac{\partial f(x_1,\cdots,x_n)}{\partial
x_i}$ is not zero polynomial. This is because if there exists
$k(0<k\leq n)$ such that $\frac{\partial
f(x_1,x_2,\cdots,x_n)}{\partial x_k}\equiv 0$, then
$$f(x_1,\cdots,x_{k-1},x_k,x_{k+1},\cdots,x_n)\equiv g(x_1,\cdots,
x_{k-1},x_{k+1},\cdots,x_n).$$ So, we can consider polynomial
$g(x_1,\cdots,x_{k-1},x_{k+1},\cdots,x_n)$. Hence with probability
1, we can get a point $(x_{1,0},x_{2,0},\cdots,x_{n,0})$ that
satisfies equation(\ref{equ:3}).
\subsection{Sampled points}

Let $X=(x_1,x_2,\cdots,x_n)$ and
$P_0=(x_{1,0},x_{2,0},\cdots,x_{n,0})$. Since $\frac{\partial
f(x_{1,0},\cdots,x_{n,0})}{\partial x_i}\neq 0$ for $i=1,\cdots,
n-1$, we can calculate $h_1>0,h_2>0,\cdots,h_{n-1}>0$ from Mean
value theorem such that
$$(sign(\frac{\partial f}{\partial x_1}(P_0)), \cdots,
sign(\frac{\partial f}{\partial
x_{n-1}}(P_0)))=(sign(\frac{\partial f}{\partial x_1}(X)), \cdots,
sign(\frac{\partial f}{\partial x_{n-1}}(X)))$$ for $X\in U(P_0)$,
where $sign()$ is sign function and let
$U(P_0)=[[x_{1,0}-h_1,x_{1,0}+h_1],\cdots,[x_{n-1,0}-h_{n-1},x_{n-1,0}+h_{n-1}]]$.

$d_i$ is denoted by the degree of polynomial $f(x_1,\cdots,x_n)$
with respect to $x_i$ for $i=1,\cdots,n-1$. Choosing $d_i+1$
distinct points in interval $[x_{i,0}-h_i,x_{i,0}+h_i]$ for
$i=1,2,\cdots,n-1$. So we get $n-1$ dimension vectors denoted by
$\Phi_{n-1}.$ And then substituting each vector $v_i\in
\Phi_{n-1}$ into $f(x_1,x_2,\cdots,x_n)$ and getting univariate
polynomial $g(x_n)$, Solving $g(x_n)=0$ within control error
$\varepsilon$ in algorithm \ref{algor:add1} and choosing the
solution $x_{n,i}$ which is the closest to $x_{n,0}$ yields a
$n$-dimension vector. Hence we get $\prod_{i=1}^{n-1}(d_i+1)$
$n$-dimension vectors denoted by $\Phi_n$.

\subsection{Getting a factor by sampled points}

In this subsection, we discuss how to get exact factors of a
polynomial. This procedure runs as follows: First a candidate set
monomials in the support of the factor is selected. If we know the
pattern of monomials of the factor, we choose a restricted set of
monomials. If we do not know the pattern, then for
$m=1,\cdots,d-1$, where $d$ is total degree of
$f(x_1,x_2,\cdots,x_n)$, we use the complete set of monomials of
degree less than or equal to $m$, denoted by $M_m$.  Next, for
every set of monomials $M_m$($m=1,2,\cdots,d-1$), some of the
above sampled points are selected. If the sampled points is less
than we want, we can refine $U(P_0)$ and get enough sampled
points. The selected sampled points must satisfy the condition of
theorem \ref{theo:1}. From the selected sampled points, formula
(\ref{equ:add1}) is used to obtain a polynomial
$g(x_1,x_2,\cdots,x_n)$ which monomial set is $M_m$. and
then transform  $g(x_1,x_2,\cdots,x_n)$ into a monic polynomial,
still denoted by $g(x_1,x_2,\cdots,x_n)$. Third, we deal with
$g(x_1,x_2,\cdots,x_n)$ in two cases: \\
Case 1: $g(x_1,x_2,\cdots,x_n)$ is not a real polynomial but a
complex polynomial. Let the selected sampled points be
$p_0,p_1,\cdots,$ from which $g(x_1,\cdots,x_n)$ is constructed.
So, their complex conjugate points $\bar p_0,\bar p_1,\cdots$, are
on the other variety $V(f_j)$. We construct a monic polynomial
$\bar g(x_1,x_2,\cdots,x_n)$ from $\bar p_0,\bar p_1,\cdots$. So
$g\bar g$ is a monic real polynomial. We deal with $g\bar g$ just
as in case 2.\\ Case 2: $g(x_1,\cdots,x_n)$ is a monic real
polynomial. We use the polynomial division proposed in
\cite{corless_B} to get $h(x_1,x_2,\cdots,x_n)$ which minimizes
$\|f-gh\|$. Let $r=\|f-gh\|$. If $r$ is large, we should add more
sampled points and extend the monomial set $M_m$ to $M_{m+1}$, and
then use formula (\ref{equ:add1}) to get a polynomial with higher
degree. If $r$ is very small, then we use algorithm \ref{alg:2} to
transform $g(x_1,\cdots,x_n)$,$h(x_1,\cdots,x_n)$ to a rational
polynomial $g'(x_1,\cdots,x_n)$ and $h'(x_1,\cdots,x_n)$
respectively. Compute $r'=\|f-g'h'\|$. If $r'=0$ then
$g'(x_1,\cdots,x_n)$ is a rational factor of $f(x_1,\cdots,x_n)$.
We continue to factor polynomial $h'(x_1,\cdots,x_n)$ on the other
variety $V(f_j)$ . If $r'\neq 0$, then from theorem
\ref{theo:practice_alg}, $g(x_1,x_2,\cdots,x_n)$ is not an
approximate rational factor but an approximate real factor of
polynomial $f(x_1,x_2,\cdots,x_n)$. We continue to factor
polynomial $h(x_1,\cdots,x_n)$ on the other $V(f_j)$. Finally, we
get all factors of $f(x_1,\cdots,x_n)$. Let
$f(x_1,\cdots,x_n)=g_1g_2\cdots g_kg_{k+1}\cdots$, and
$\Phi=\{g_1,g_2,\cdots,g_k\}$ are not rational factors, the others
are rational factors. We compute products of two distinct
polynomials of $\Phi$: $g_{ij}=g_i*g_j (i\neq j, and\,\, i,j\leq
k)$. Check if every $g_{ij}$ is an approximate rational factor,
Whenever it is so, keep $g_{ij}$ as a rational factor and remove
$g_i,g_j$ from $\Phi$; After finishing to deal with products of
two distinct polynomials in $\Phi$, we compute products of three
distinct polynomials of $\Phi$:\ $g_{ijm}=g_ig_jg_m (i\neq j,j\neq
m ,m\neq i)$, and check if every $g_{ijm}$ is an approximate
rational factor. Whenever $g_{ijm}$ is so, keep it as a rational
factor and remove $g_i,g_j,g_m$ from $\Phi$; and so on, until all
rational factors of $f(x_1,x_2,\cdots,x_n)$ are obtained.

\section{Experimental results}

The following five polynomials are randomly generated by
\emph{randpoly} command in \emph{Maple}. Our algorithm is
implemented in \emph{Maple}. Compared our algorithm with
\emph{factor} command in \emph{Maple} in the platform of Maple 9.5
and PIII 1.0G, 256M RAM, The running time of the five examples are
as follows:

\textbf{Example 1.} This polynomial is with four variables, 46376
terms, and of degree 30. The 15 factors are of degree 2. The
\emph{factor} command in \emph{Maple} costs 430.108 seconds, our
algorithm costs 278.375 seconds.

\textbf{Example 2.} The polynomial is with four variables, 52360
terms, and of degree  31. The 21 factors are of degree 1 or 2. The
\emph{factor} command in \emph{Maple} costs 667.438 seconds to
factor the polynomial, our algorithm 307.922 seconds.

\textbf{Example 3.} The polynomial is with four variables, 52360
terms, and of degree 32. The 32 factors are of degree 1. The
\emph{factor} command in \emph{Maple} computed 7200 seconds and
gave no result, our algorithm  250.436 seconds.

\textbf{Example 4.} The polynomial is with three variables, 17296
terms, and of degree  45. The 15 factors are of degree 3. The
\emph{factor} command in \emph{Maple} use 261.468 seconds to
factor the polynomial, our algorithm 173.265 seconds.

\textbf{Example 5.} The polynomial is with three variables, 37820
terms, and of degree  59. The 25 factors are of degree 2 or 3. The
\emph{factor} command in \emph{Maple} used 1683.655 seconds to
factor the polynomial, our algorithm used 710.045 seconds.
\section{Conclusion}
\begin{itemize}

\item Our algorithm first need a initial point. We just get $n-1$
floating-point number at random, and then substitute them into a
polynomial to be factored and obtain a univariate polynomial.
Solving the univariate polynomial for the last variable $x_n$
yielding $d_n$ solutions. In the above discussion, it seems only
to take one solution and throw away the other solution. In fact,
we should keep these $d_n-1$ points as initial points on the other
varieties. For the same reason, we should keep other $d_n-1$
solution whenever get a sampled point in neighborhood $U_0$. Once
we get sampled points in neighborhood $U_0$, we obtain $d_n-1$
sets of sampled points on the other $d_n-1$ varieties of the
original polynomial.

\item Our algorithm is to get exact factorization of polynomials
by approximate method, so its efficiency is higher than symbolic
factorization method when a polynomial is with more variables and
higher degree. However its efficiency is lower than symbolic
factorization when a polynomial is with less variables and lower
degree. In order to take advantage of both numerical and symbolic
factorization algorithm, we can factor a polynomial with more
variables and higher degree as follows: First, we use our
algorithm to get some factors of the polynomial and remove these
factors from the polynomial. When the remaining polynomial is with
less variables and lower degree, then, we use symbolic
factorization algorithm to get the remaining factors. In fact, Our
algorithm is implemented by this idea.

\item Noting that interpolating formula (\ref{equ:add1}) is lower
efficient when the factors of the original are with more variables
and higher degree. So we should improve interpolating formula
(\ref{equ:add1}) further.
\end{itemize}

\end{document}